\documentclass[10pt,a4paper]{article}
\usepackage{graphicx}
\usepackage[unicode=true]{hyperref}
\usepackage{amsmath}
\usepackage{hyperref}
\hypersetup{
    colorlinks=true,
    linkcolor=blue,
    filecolor=magenta,      
    urlcolor=cyan,
}
\usepackage{amssymb}
\usepackage{natbib}
\usepackage[frozencache]{minted}
\usepackage[author=IB]{pdfcomment}
 \def\style{newapa} 
\begin{document}

\title{Conceptual Mathematics via Literate Programming}
\author{Ian Benson, Jim Darby, Neil MacDonald and Jesse Sigal}
\maketitle

Recent developments in computer programming and in mathematics suggest that there is a strong case for a new way of introducing programming to enhance the learning of school mathematics. We believe that a programming approach based on both types and functions would make a vastly improved contribution to learning  mathematics than the less successful use of conventional computer programming in Scratch \citep{Boylan:ly,Rodriguez-Martinez:2020to}. 

A key to this symbiosis between school mathematics and informatics is Knuth's literate programming approach to coding. The main idea is to let the program grow in natural stages, with its parts presented in roughly the order that they might have been created by a programmer who isn't especially clairvoyant. Knuth wanted to develop a language for ``top-down'' programming, where a top-level description of the problem situation is given first and successively refined. \textit{What he found was that a program is best thought of as a web instead of a tree}. A hierarchical structure is present, but the main thing about a program is its structural relationships. ``A complex piece of software consists of simple parts and simple relations between those parts; the programmer’s task is to state those parts and those relationships, in whatever order is best for human comprehension— not in some rigidly determined order like top-down or bottom-up'' \citep[p. 11]{Knuth:1960vb}. 

A key goal of our application of literate programming is to enable the reader to validate by inspection (and inline execution) that the code meets expectations. A conceptual model sets out these design objectives. It will often be recorded in symbols, notation and/or structured drawings. A clear separation of the theoretical or conceptual model from its (numerical) specification and code helps the designer, and her users, to validate how the conceptual model is implemented. That is, to determine that the implementation satisfies the human stakeholders’ requirements, while also verifying that the specification of the program is met by the implementation.

In this article we describe a case study of a comparative programming exercise, in which solutions to a Josephus problem are developed in Python and Haskell by teachers of mathematics and computer science, and program behaviour is modelled using conceptual mathematics. We show how number theoretic mathematics (modulo arithmetic) contributes to an imperative language solution in Python, while  conceptual mathematics helps to record the isomorphic relationship between the dynamic behaviour of specification and implementation.

\section{Conceptual mathematics}

Conceptual mathematics, or Category theory, is becoming a central feature of modern pure mathematics. It gets its power from the ability to organize and refine abstractions, to find common aspects between many different mathematical structures and to facilitate communication between different communities of mathematicians and between mathematicians, computer scientists, logicians, and physicists \citep{Fong:2019di}. William Lawvere and Stephen Schanuel have written an accessible introduction for secondary school mathematics teachers. They argue that ``Everyone who wants to follow the applications of mathematics to twenty-first century science should know (these) ideas and techniques.'' \citep{Lawvere:2004zk}.

\section{Josephus's Problem}

Mathematics teachers at an MEI event worked on the following problem:
\begin{quote}
When dealing to the Roman Army, the term decimate meant that the entire unit would be broken up into groups of ten soldiers, and lots would be drawn. The person who was unlucky enough to draw the short straw would be executed by the other nine members of his group. The bloodthirsty Roman Centurion Carnage Maximus decided to apply this to his prisoners, with a few gruesome differences. Rather than kill every tenth prisoner and allow the rest to live, he is going to leave only one prisoner alive and kill all of the others. There are 100 prisoners and they stand in a circle. He goes round killing every 10th prisoner until one is left. Which one?
\end{quote}


This is a variant of an ancient problem named for Flavius Josephus, a famous historian of the first century. Legend has it that he wouldn’t have lived to become famous without his mathematical talents. During the Jewish-Roman war, he was among a band of 41 Jewish rebels trapped in the \href{https://en.wikipedia.org/wiki/Josephus_problem}{Siege of Yodfat} by the Romans.  Preferring to die rather than to be captured, the rebels decided to form a circle and, proceeding around it, to kill every third remaining person until no one was left. But Josephus, along with an unnamed co-conspirator, wanted none of this nonsense about suicide;  so he quickly calculated where he and his friend should stand in the vicious circle \citep[p. 8]{Ronald-L.-Graham:1993kl}. 

Knuth estimated that this is a problem of moderate difficulty and/or complexity which could take two hours to solve satisfactorily, or even more if the TV is on \citep[p. 162]{Knuth:1997km}. The reader is encouraged to take a break to consider how she might attempt to solve the problem by hand, by programming in a pseudo-code such as AQA or in your favourite language in the general case of $n$ men in a circle with every $m’th$ man executed in turn \citep{AQA:2jk}.

\subsection{A `WEB' program primer}

Knuth illustrated his approach by publishing the code for his TEX and METAFONT Pascal programs in WEB. Our first algorithm was developed in Python by Michael Jones an NCCE computer science champion and revised below as annotated with the `WEB' formalism. An extended version of this example with a complete set of cross-referencing notes will be found at ...

\textbf{Note 1.} We introduce \mintinline{python}{prisoners} and \mintinline{python}{index} variables with the following life cycles:\\
$\langle$ Procedures for data manipulation 1 $\rangle$ $\equiv$
\begin{table}[htbp]
\begin{tabular}{p{10.0ex}p{55.0ex}}  \hline
\mintinline{python}{prisoners} & create a new \mintinline{python}{prisoners} list (line 9)\\
& kill the prisoner at \texttt{prisoners[index]} and remove a dead prisoner (line 7)\\
& quit when only one prisoner remains (lines 5)\\ \hline
\mintinline{python}{index} &  create a new \mintinline{python}{index} (line 3)\\
& \mintinline{python}{index = (pos + index) 
\end{tabular}
\end{table}

A python version of this code is implemented in lines 2 through 7.

\textbf{Note 2.} The program represents the living prisoners as a ordered list of integers [1..100], with a cursor (\mintinline{python}{index}) pointing to the next prisoner to kill. \mintinline{python}{len(prisoners)} is the length of the current list of live prisoners. It changes after each kill. \mintinline{python}{index} is recalculated using modulo  \mintinline{python}{len(prisoners)} arithmetic. The program exits when only one prisoner remains \\
$\langle$ The main program 2  $\rangle$ $\equiv$
\begin{table}[htbp]
\begin{tabular}{p{0.1ex}p{75.0ex}} 
 & initialise \mintinline{python}{prisoners}\\
 & call local procedure \mintinline{python}{removeTen}\\
 & quit when \mintinline{python}{removeTen} returns with only single prisoner remaining (line 5) \\
 \end{tabular}
\end{table}

A python version of this code is implemented in lines 9 through 11.

\subsection{Python implementation, \texttt{romans.py}}

\begin{minted}[linenos]{python}
def removeTen(prisoners):
    pos = 10 - 1
    index = 0

    while len(prisoners) > 1:
        index = (pos + index) % len(prisoners)
        prisoners.pop(index)

prisoners = list(range(1,101))
removeTen(prisoners)
print(prisoners)
\end{minted}

\section{An experiment in collaborative programming}

\begin{figure*}
  \includegraphics[width=.7\textwidth]{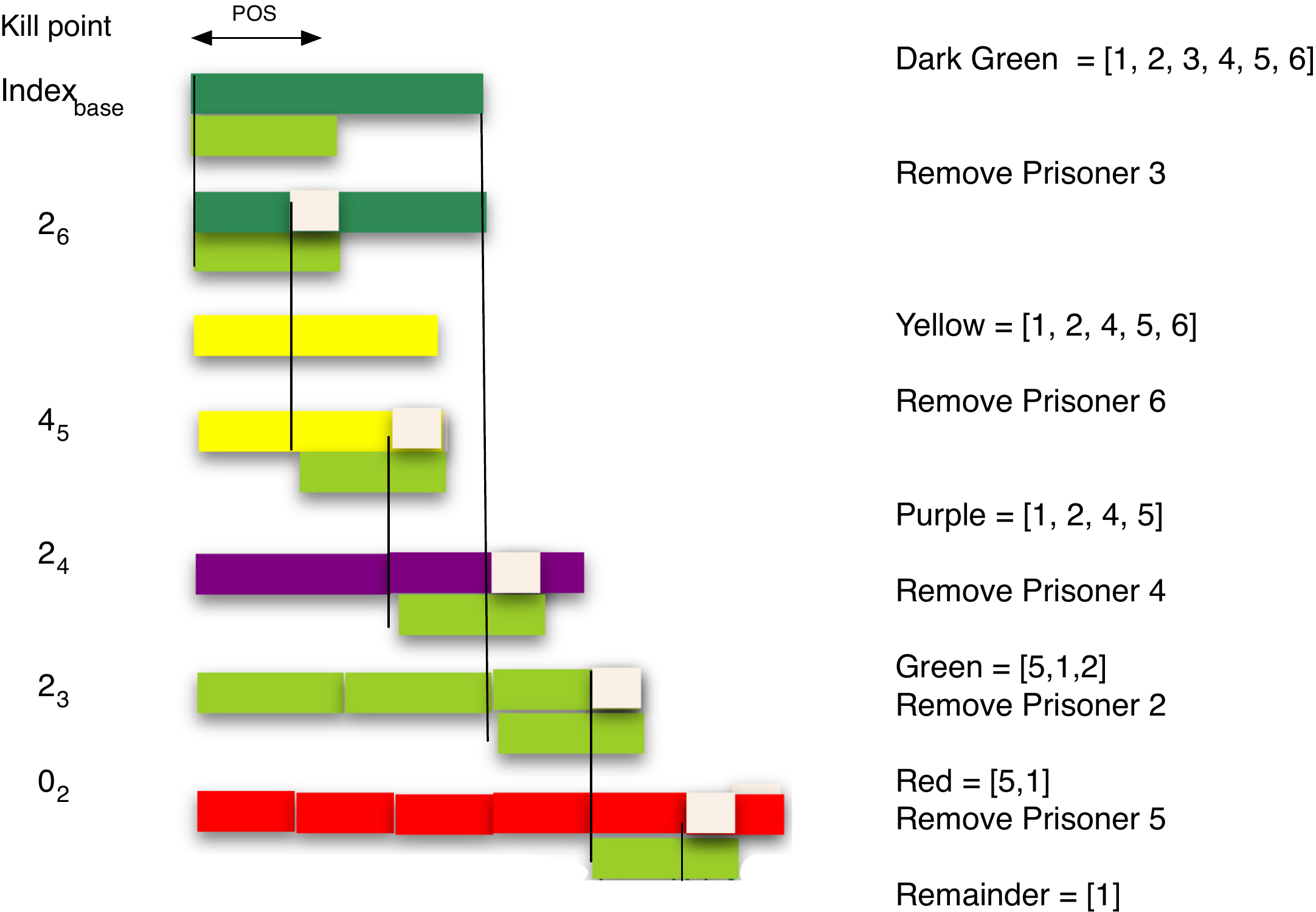}
\caption{A circle of 6 prisoners, killing every 3rd.}
\label{fig:CUI1}       
\end{figure*}

An initial version of the solution, which records the mathematics teachers intuition, is shown in Figure~\ref{fig:CUI1}. This illustrates the behaviour of a model with a circle of 6 prisoners with every 3rd killed. Mathematics teachers given this exercise developed an intuition that the solution might have something to do with number systems in different bases. A computer science teacher implemented a similar idea in \texttt{romans.py}. The intuition concerning number bases is limited. For example, when we start with a large number of prisoners, say 20, we are not working in base 20 in the sense that there are 20 unique symbols for a new one's place. Each kill cycle is shown as a train of single coloured rods. The train is recalibrated with each step, measuring it with the predecessor rod. The kill position is also recalibrated relative to the last man killed. The mapping between the rod name and the list of remaining prisoners is shown on the right of the figure. It is left as an exercise for the reader to account for the change in the cycling of the prisoners in the step from purple to green. Several other implementations of the solution were developed and shared though github and repl.it.

\section{Haskell solution \texttt{romans.hs}}
\label{romans2}

\subsection{Specification Model: Function type signatures}

A Haskell solution to the problem is split into two main portions.
We first create a conceptual model of a non-empty circle of prisoners via the abstract data type \mintinline{haskell}{CircleOf a}.
A instance of \mintinline{haskell}{CircleOf a} represents a circular arrangement of values of type \texttt{a} with one element as the current focus.
We now specify a collection of basic operations to work with \mintinline{haskell}{CircleOf a}.
The function \mintinline{haskell}{mkCircleOf} takes an value for the focus and a list of values to be after it and creates a \mintinline{haskell}{CircleOf a}.
Additional functions defined are: \mintinline{haskell}{current} which returns the element currently in focus, \mintinline{haskell}{isSingleton} which returns \mintinline{haskell}{True} if the circle contains exactly one element, \mintinline{haskell}{next} which moves the focus to the next element, and \mintinline{haskell}{remove} which removes the element at the focus while placing the next element in focus.

The model enables us to write a function \mintinline{haskell}{removeNth} which takes an integer $n$ and a circle of elements which calls \mintinline{haskell}{next} $n-1$ times and then calls \mintinline{haskell}{remove}, removing the $n^\text{th}$ element from the original focus.
Finally, we can write the main function \mintinline{haskell}{romans} which takes the number of prisoners and which $n^\text{th}$ prisoner to kill and returns it.
We use the library function \mintinline{haskell}{until} which takes three arguments, where \mintinline{haskell}{until p f x} applies \texttt{f} (possibly zero times) to \texttt{x} until \texttt{p} is true. 

\begin{minted}[linenos]{haskell}
mkCircleOf :: (a, [a]) -> CircleOf acircle 
current :: CircleOf a -> a
isSingleton :: CircleOf a -> Bool
next :: CircleOf a -> CircleOf a
remove :: CircleOf a -> CircleOf a

removeNth :: Int -> CircleOf a -> CircleOf a
removeNth 1 circle = remove circle
removeNth n circle = removeNth (n-1) (next circle)

romans :: Int -> Int -> Int
romans numPrisoners n =
  let prisoners = mkCircleOf (1, [2..numPrisoners])
  in current (until isSingleton (removeNth n) prisoners)

main = print (romans 100 10)
\end{minted}

\subsection{Implementation Model: Defining Equations}

The conceptual model described is implemented below.
The data type \mintinline{haskell}{CircleOf a} is a pair of the focussed element and a list of the other elements in order.

\begin{minted}[linenos,firstnumber=17]{haskell}
data CircleOf a = C a [a] deriving Show
mkCircleOf (x, xs) = C x xs
current (C x _) = x
isSingleton (C _ []) = True
isSingleton _        = False
next (C x []    ) = C x []
next (C x (y:ys)) = C y (reverse (x : reverse ys))
remove (C x []    ) = C x []
remove (C x (y:ys)) = C y ys
\end{minted}

\subsection{An application of conceptual mathematics} 

\begin{figure*}
\centering
\includegraphics[width=.9\textwidth]{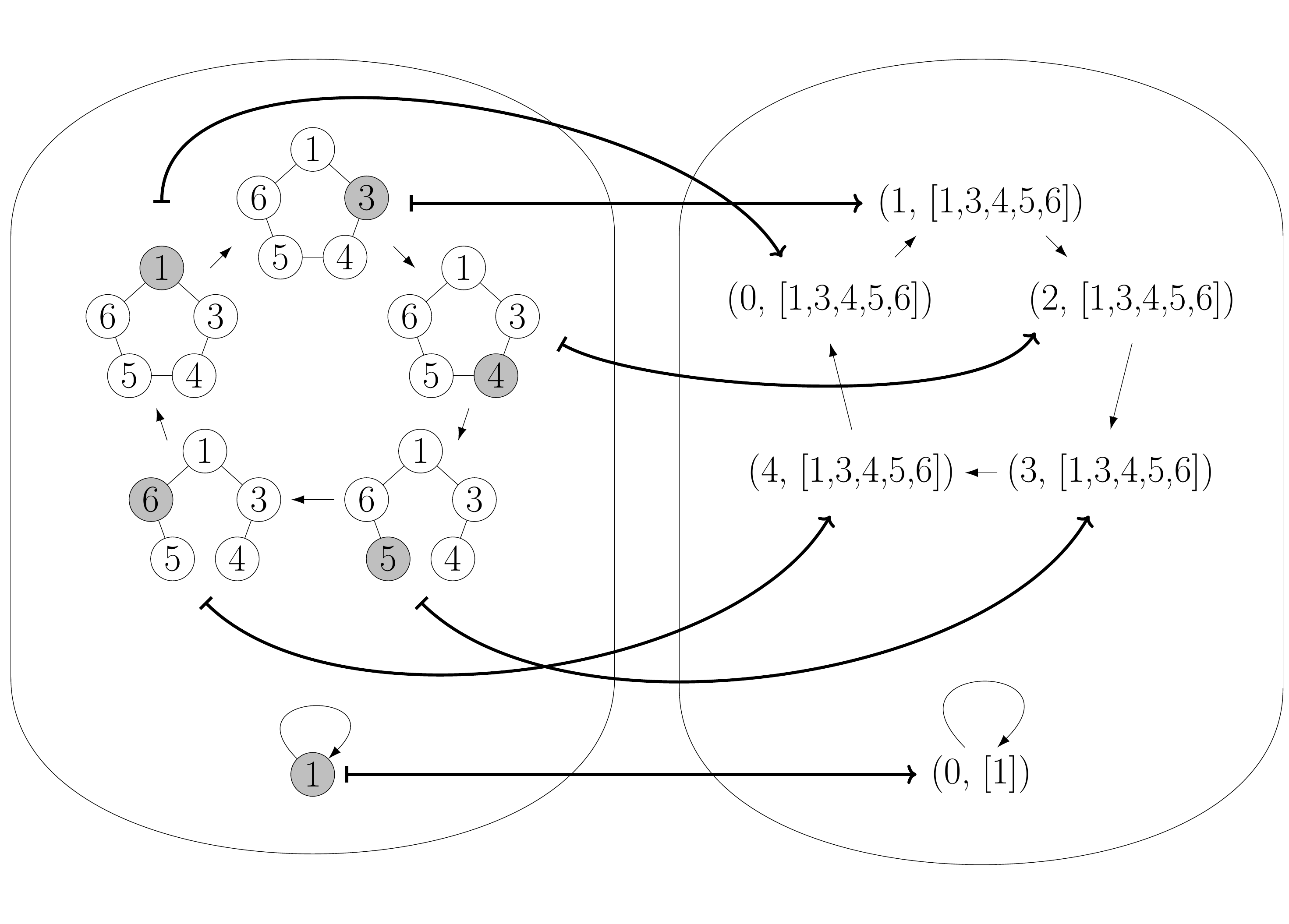}
\caption{An internal diagram of an isomorphism between program behaviours}
\label{Morphism}       
\end{figure*}

For a set $A$, a map $\alpha \colon A \to A$ from $A$ to itself is an \emph{endomap of $A$}.
A set $A$ paired with an endomap $\alpha$, written as $(A, \alpha)$, can be seen as a discrete dynamical system with $A$ giving the states and $\alpha$ giving the one-step dynamics, i.e. a state $x \in A$ transitions to $\alpha(x)$ after one step.
Let $(A, \alpha)$ and $(B, \beta)$ be discrete dynamical systems, then a map $f \colon (A, \alpha) \to (B, \beta)$ between them is a map $f \colon A \to B$ such that $f \circ \alpha = \beta \circ f$. 
\cite[p. 152]{Lawvere:2004zk} propose an internal diagram notation to show maps (also called morphisms) between discrete dynamical systems.
For a dynamical system $(A, \alpha)$, we draw some elements $x \in A$ and a line from $x$ to $\alpha(x)$. 
Let $H$ denote the set of instances of \mintinline{haskell}{CircleOf Int} with
distinct elements between $1$ and $6$.
We can then view $(H, \mathtt{next})$ as a discrete dynamical system.
Let $P$ denote the set of pairs of a number $0$ to $5$ and a list of length at most $6$ with distinct elements which are integers $1$ to $6$.
The set $P$ represents possible values of the pair $(\mathtt{index}, \mathtt{prisoners})$ in the Python solution.
Let $\gamma$ denote the action of line $6$ of the Python program.
We can then view $(P, \gamma)$ as a discrete dynamical system.
Figure~\ref{Morphism} illustrates a morphism $f \colon (H, \mathtt{next}) \to (P, \gamma)$ on example elements.
In fact, if we invert the mapping, we obtain another map of dynamical systems, and so $(H, \mathtt{next})$ and $(P, \gamma)$ are isomorphic as dynamical systems. 
Thus, the conceptual model of discrete dynamical systems allows us to compare programs across different languages.

\section{Conclusion}

Francis Su cites Bill Thurston as saying that ``the real satisfaction from mathematics is in learning from others and sharing with others. All of us have clear understanding of a few things and murky concepts of many more'' \citep{Su:2021az}.  In his book Su offers an extended example ``to illustrate many of the powers of mathematics — interpretation, definition, quantification, abstraction, visualization, imagination, creation, strategization, modeling, multiple representations, generalization, and structure identification. Anyone who learns mathematics will grow skilled in these powers. These are virtues that enable the creative power of `stuff making and sense making' ''. 

In this case study we have shown Su’s `mathematical powers’ at work through literate programming:

\begin{description}
  \item[definition:] abstract data type definition of \mintinline{haskell}{CircleOf} (\mintinline{haskell}{mkCircleOf}, \mintinline{haskell}{current}, \mintinline{haskell}{isSingleton}, \mintinline{haskell}{next}, \mintinline{haskell}{remove})
  \item[abstraction:] the data type \mintinline{haskell}{CircleOf}, the functions \mintinline{haskell}{romans} and \mintinline{haskell}{removeNth}, the data type choice for \mintinline{python}{prisoners}
  \item[visualisation:] figures \ref{fig:CUI1} and \ref{Morphism}
  \item[creation:] conceptual modelling
  \item[strategization:] for the Haskell solution, the use of \mintinline{haskell}{until}, for the Python solution, the function \mintinline{python}{removeTen}
  \item[modelling:] lists in Python and \mintinline{haskell}{CircleOf a} in Haskell
  \item[generalisation:] generic functions in Python (\mintinline{python}{pop(index)}) and Haskell (\mintinline{haskell}{until})
  \item[structure identification:] each main program
\end{description}

\subsection{Postscript}

To validate this case study we shared an earlier version of this article with mathematics teachers who were familiar with using versions of Josephus’s problem with Y7 and  Y12 pupils (on the white board). We also researched the literature on mathematical and programming solutions to the problem. What we discovered was that the problem was accessible to learners across secondary school and could be used to encourage the exploration of both mathematical and computing solutions. Pupils could even explain their reasoning inductively in the special case of $n$ a power of 2, and $m=2$ (In this case the wording of the problem was changed to the less gruesome `musical chairs' terminology of in (live) out (dead)). We discovered that computer scientists have worked extensively on designing more efficient algorithms than the ones we discuss here. Knuth himself with colleagues has developed ingenious solutions to generalizations of the Josephus problem \citep{Knuth:1997dl, Ronald-L.-Graham:1993kl,Woodhouse:1978nq,Lloyd:270sh}.

\section{About the authors}
Ian Benson, Jim Darby, Neil MacDonald and Jesse Sigal are members of the ATM Functional Programming and Computer Algebra Working Group. They are grateful for comments by Jim Thorpe, David Vaccaro and Richard Hunt on earlier drafts of this article.  This pre-print is under review for publication. The corresponding author is \href{mailto:Ian.Benson@Roehampton.ac.uk}{Ian.Benson@Roehampton.ac.uk}.

\bibliography{NewBib.bib}   
\bibliographystyle{\style}

\end{document}